\documentclass[12pt]{article}

\usepackage{amssymb}
\usepackage{amscd}
\usepackage{amsmath}

\newtheorem{theorem}{Theorem}
\newtheorem{proposition}[theorem]{Proposition}

\newtheorem{definition}[theorem]{Definition}
\newtheorem{exam}[theorem]{Example}

\newtheorem{rem}[theorem]{Remark}

\numberwithin{theorem}{section}

\renewcommand{\subsection}[1]{\vspace{5 mm}
\noindent \bf{\Large{#1}} \rm
\par}

\newcommand{\pp}{\mathbb{P}}

\newcommand{\cc}{\mathbb{C}}
\newcommand{\rr}{\mathbb{R}}
\newcommand{\zz}{\mathbb{Z}}

\newcommand{\lik}{\mathfrak{k}}
\newcommand{\lig}{\mathfrak{g}}
\newcommand{\lih}{\mathfrak{h}}

\newcommand{\lit}{\mathfrak{t}}

\newcommand{\git}{/\!\!/}

\newcommand{\Cal}{\mathcal}

\begin{document}

\title{Intersection cohomology of representation spaces
of surface groups}
\author{Young-Hoon Kiem\\
Department of Mathematics\\ 
Stanford University\\
Stanford, CA 94305\\
kiem@math.yale.edu}
\date{}
\maketitle

\section{Introduction}
Let $\Sigma$ be a closed Riemann surface of genus $g\ge 2$
and $\pi$ be its fundamental group.
The representation space $$X(G)=\mathrm{Hom}(\pi,G)/G$$
of $\pi$ into a compact
connected Lie group $G$ has been playing important roles in
differential geometry (the moduli space of flat connections),
algebraic geometry (the moduli space of holomorphic principal
bundles) and topology (Casson's invariant).
This is a compact pseudomanifold which can be thought of
as the symplectic reduction of a Hamiltonian $G$-space,
the extended moduli space. The purpose of this paper is to
study the local structure of $X(G)$ in terms of the symplectic
geometry of the extended moduli space and to show 
by the splitting theorem \cite{kiem, kw} that its intersection cohomology
is embedded into the $G$-equivariant cohomology of 
$\mathrm{Hom}(\pi,G)$, which can be computed by the
equivariant Morse theory of Atiyah-Bott \cite{ab2} \S10
at least in principle. The image is identified with
a subspace $V^*_U$ truncated locally and 
this enables us to compute the intersection cohomology
as a graded vector space with a nondegenerate pairing.

The extended moduli space $\mathcal{N}^{\lig}$ for a simply
connected compact Lie group $G$ was defined and studied by
L. Jeffrey \cite{j}. Roughly, it is a ``partial
reduction'' of the space of connections by the based
gauge group and thus carries a symplectic structure
on the smooth part together with a (residual) $G$ action.
$X(G)$ is now the symplectic reduction of this partial
reduction. In \cite{JK2}, L. Jeffrey and F. Kirwan applied
their nonabelian localization principle to the extended
moduli space and proved Witten's formulas for the intersection
numbers of the moduli spaces of vector bundles of coprime rank
and degree. In \S2, we review the definition and some key
results for $\mathcal{N}^{\lig}$. Then we determine the
symplectic slice at a point in the zero set of the moment map.
As is well-known, at a smooth point in $X(G)$, the tangent 
space is isomorphic to the symplectic slice. Our description
of the symplectic slice is consistent with W. Goldman's 
description of the tangent space \cite{Gold}.

Given a singular reduction, R. Sjamaar and E. Lerman provided a useful local description \cite{SL} and showed that a singular reduction
is a stratified space, symplectic on each stratum. In \S3,
we apply their results to study the local structures of $X(G)$.
Locally, $X(G)$ looks like the symplectic reduction of 
the symplectic slice by a linear action of the stabilizer.
This gives us a local description of $X(G)$, consistent
with the result of Goldman-Millson in \cite{GolMil},
where they obtained analytic results through deformation theory.

Since $X(G)$ is a singular pseudomanifold, intersection cohomology
is an interesting invariant, at least as important as the ordinary
cohomology because of its nice properties including Poincar\'e duality
in spite of singularities. Generalizing the results of \cite{kiem},
we showed in \cite{kw}
that the intersection cohomology of a reduced space $M\git G$
of a Hamiltonian $G$ space $M$ with proper moment map $\mu$
is naturally isomorphic to a subspace, called the truncated
equivariant cohomology, of the equivariant cohomology
$H^*_G(\mu^{-1}(0))$ under an assumption named 
\emph{weakly balanced}. This condition was devised to make
the codimensions of the unstable strata in the blow-ups 
in the partial desingularization process, large enough. 
The truncated equivariant cohomology is defined by truncating
the equivariant cohomology along each stratum and the splitting
theorem\footnote{The isomorphism is a splitting of 
the Kirwan map $\kappa:H^*_G(\mu^{-1}(0))\to IH^*(\mu^{-1}(0)/G)$.}
enables us to compute the Betti numbers and intersection pairing
of $IH^*(M\git G)$ in terms of the equivariant cohomology.
In \S4, we recall this theorem and show that the extended
moduli space satisfies the assumption. Therefore, we can
compute the intersection cohomology of $X(G)$ in terms of 
the $G$ equivariant cohomology of $\mathrm{Hom}(\pi,G)$.

The case where $G=SU(2)$ is discussed in \S5. We can compute
the intersection cohomology groups of $X(SU(2))$ by using
the structure theorem of $H^*_{SU(2)}(\mathrm{Hom}(\pi,SU(2)))$
from \cite{kiem4}.

If the Lie group $G$ is not simply connected, $X(G)$ is not
connected. In fact, $\pi_0(X(G))\cong \pi_1(G)$. Let $\tilde{G}$
denote the universal cover of $G$ with (finite) fiber $C$. 
Then there is a covering space
$X(\tilde{G})\to X(G)$ onto a component we denote by $X(G)_0$.
Then the intersection cohomology of $X(G)_0$ is the $C^{2g}$
invariant part of the intersection cohomology of $X(\tilde{G})$
which can be computed by the splitting theorem.

Every cohomology group
in this paper has complex coefficients and 
every intersection cohomology has middle perversity.

\textbf{Acknowledgements.} I am grateful to Professors Ronnie Lee, Lisa
Jeffrey, Frances Kirwan and Jon Woolf for useful discussions.

\section{Extended moduli space}
Let $G$ be a simply connected compact Lie group and $\Sigma$ be
a closed Riemann surface of genus $g\ge 2$. Let $\pi$
denote the fundamental group of $\Sigma$. The representation space
$$X(G):=\mathrm{Hom}(\pi,G)/G$$
can be realized as the symplectic reduction of a Hamiltonian $G$
space $\mathcal{N}^{\lig}$, called the \emph{extended moduli space}. 
This was defined and studied by L. Jeffrey in \cite{j}.

\begin{definition} Let $\lig=Lie(G)$. We define
$$\mathcal{N}^{\lig}:=\{(a_1,\cdots,a_{2g},\Lambda)\in
G^{2g}\times\lig | \, \prod[a_i,a_{i+g}]=e^{\Lambda}\}$$
More generally, for each central element $c$ in $G$, we define
$$\mathcal{N}^{\lig}_c:=\{(a_1,\cdots,a_{2g},\Lambda)\in
G^{2g}\times\lig | \, \prod[a_i,a_{i+g}]=ce^{\Lambda}\}.$$
\end{definition}

Let $P:G^{2g}\times \lig\to G$ be the map defined by 
$$P(a_1,\cdots,a_{2g},\Lambda)=\prod [a_i,a_{i+g}]e^{-\Lambda}.$$
Then $\mathcal{N}^{\lig}=P^{-1}(1)$. 
$dP_{(a,\Lambda)}$ is surjective if $P(a,\Lambda)=1$ and 
$\Lambda=0$ because $d(exp)_0$ is surjective. Hence,
$\mathcal{N}^{\lig}$ is smooth in a neighborhood of 
$$\mathrm{Hom}(\pi,G)=\{(a_1,\cdots,a_{2g})\in G^{2g}\,|\,
\prod [a_i,a_{i+g}]=1\}$$
It is not hard to see that $\mathcal{N}^{\lig}$ is connected.
Also, from the work of A. Ramanathan \cite{Ram} and Atiyah-Bott \cite{ab2},
it is easy to see that $\mathrm{Hom}(\pi,G)$  is connected.

The extended moduli space can be considered as a moduli space
of flat connections as follows. Delete an open disk $D$ with 
boundary $S$ from $\Sigma$ and denote the resulting surface by 
$\Sigma'$. Fix a parametrization of $S$ by a parameter $s\in \rr/\zz$
and a diffeomorphism of a neighborhood of $S$ to $\rr/\zz\times
[0,\epsilon)$. Let $\mathcal{A}_{flat}^{\lig}(\Sigma')$ denote
the set of flat connections on $\Sigma'$ whose restriction to a 
neighborhood of $S$ is $\alpha\, ds$ for some $\alpha\in\lig$.
Let $\mathcal{G}_0(\Sigma')$ be the space of maps $g:\Sigma'\to G$
such that $g=I\in G$ in a neighborhood of $S$. Then, as usual,
holonomy defines a map $\mathcal{A}_{flat}^{\lig}(\Sigma')/
\mathcal{G}_0(\Sigma')\to \mathcal{N}^{\lig}$. We
quote the following theorem from \cite{j}.

\begin{theorem} $\mathcal{N}^{\lig}$ is homeomorphic to 
$\mathcal{A}^{\lig}_{flat}(\Sigma')/\mathcal{G}_0(\Sigma')$.
\end{theorem}

From this gauge theoretic  description, we see that the 
tangent space at a smooth point in $\mathcal{N}^{\lig}$ is
isomorphic to the first cohomology $\tilde{H}^1$ of the complex
(see \cite{j}) 
\begin{equation}
\begin{CD}
{\Omega^0_c(\Sigma')\otimes \lig}@>{d_A}>> {\Omega^{1,\lig}(\Sigma')}
@>>>{\Omega^2_c(\Sigma')\otimes \lig}
\end{CD}
\label{tangcx}
\end{equation}
where $\Omega^i_c$ is the space of $i$-forms with compact support
in the interior of $\Sigma'$ and $\Omega^{1,\lig}(\Sigma')=\{
b\in \Omega^1(\Sigma')\otimes \lig\,|\, b=\alpha ds
\text{ in a neighborhood of } S \text{ for some }\alpha
\in \lig\}$. Notice that  $\Omega^{1,\lig}(\Sigma')\cong
\Omega^1_c(\Sigma')\otimes \lig \oplus \lig$
by choosing a 1-form whose restriction to a neighborhood of $S$
is $ds$. Hence, the index
of the complex is $-2g (\dim G)$ since the Euler characteristic
of $\Sigma'$ is $1-2g$. Therefore, the dimension of $\mathcal{N}^{\lig}$
is $2g (\dim G)$.

Let $a,b\in\Omega^{1,\lig}(\Sigma')$ and define a pairing $\omega$
on $\Omega^{1,\lig}(\Sigma')$ by 
\begin{equation}\label{omegaN}
\omega(a,b)=\int_{\Sigma'}\langle a,b\rangle\end{equation}
where $\langle,\rangle$ is the nondegenerate invariant metric on $\lig$.
By Stokes' theorem, $\omega$ defines a 2-form on $\tilde{H}^1$, which
is shown to be symplectic in \cite{j}. Thus near $\mathrm{Hom}
(\pi,G)$, $\mathcal{N}^{\lig}$ is a symplectic manifold.

For $(a_1,\cdots, a_{2g},\Lambda)\in \mathcal{N}^{\lig}$ and $h\in G$,
consider the action $$h\cdot (a_1,\cdots,a_{2g},\Lambda)=(ha_1h^{-1},\cdots,
ha_{2g}h^{-1},h\Lambda h^{-1}).$$
This makes $\mathcal{N}^{\lig}$ a $G$ space. In terms of gauge theory,
this action can be described by the corresponding action on the 
space of connections. Suppose $\zeta\in \lig$ and let 
$a\in\Omega^{1,\lig}(\Sigma')$ be a $d_A$-closed 1-form
whose restriction to a neighborhood of $S$ is $\alpha\, ds$.
Then \begin{equation}
\omega(d_A\zeta, a)=\int_{\Sigma'}\langle d_A\zeta,a\rangle =
\langle \zeta,\alpha\rangle.\label{momcom}\end{equation}

Let $\mu:\mathcal{N}^{\lig}\to \lig$ be the map defined by
$\mu(a_1,\cdots,a_{2g},\Lambda)=-\Lambda$. Then the above
computation shows that $\mu$ is a moment map\footnote{The sign
depends on the convention.} for the $G$ action
on a neighborhood $U=\mu^{-1}(B)$, for a small open ball $B$ around 0,
of $\mathrm{Hom}(\pi,G)=\mu^{-1}(0)$. In particular, the representation
space $X(G)$ is the symplectic reduction $U/\!/G$ of $U$ by the $G$ action.

Near a point $\phi=[A]$ in $\mu^{-1}(0)$, $\mathcal{N}^{\lig}$
is diffeomorphic to $\tilde{H}^1$ of (\ref{tangcx}). The orbit
direction is given by  $d_A\zeta$ for $\zeta\in \lig$ and thus 
by (\ref{momcom}) its symplectic orthogonal complement in 
$\tilde{H}^1$ is precisely $H^1_c(\Sigma';d_A)\cong
H^1(\Sigma,D;d_A)$. Notice that the local system extends
to $\Sigma$ because $\phi=[A]\in \mu^{-1}(0)$. 
 The infinitesimal action of $\lig$ is given
by the image of the map
$$H^0(D;d_A)\to H^1(\Sigma,D;d_A)$$
Therefore, from the long exact sequence, we deduce the following theorem.
\begin{theorem}
The symplectic slice at a point $\phi=[A]$ in $\mu^{-1}(0)\subset \mathcal{N}^{\lig}$
is $H^1(\Sigma;d_A)=H^1(\pi;\lig_{Ad\phi})$.
\end{theorem}

Hence, the tangent space at a smooth point $\phi$ of $X(G)$
is isomorphic to $H^1(\pi;\lig_{Ad\phi})$. This is consistent
with Goldman's result \cite{Gold}. From
(\ref{omegaN}), the symplectic structure  on
the slice is given by the cup product and $\langle,\rangle$.


\section{Local description}
In this section, we apply the results of Sjamaar-Lerman \cite{SL}
to study the local structure of $X(G)$ and $U\subset
\mathcal N^{\lig}$ with the $G$ action.

Let $M$ be a Hamiltonian $G$ space with moment map $\mu$.
We recall the following local normal form theorem.
\begin{theorem} \cite{SL}
Let  $p\in \mu^{-1}(0)$ and $\hat{W}_{p}$
be the symplectic slice of the orbit $Gp$. Then a neighborhood of 
the orbit is equivariantly symplectomorphic to a neighborhood of
the zero section of $$G\times_{Stab\,p}((\lig/\lih)^*\times\hat{W}_{p})$$
with moment map $\mu(a,\xi,v)=Ad^*(a)(\xi+\mu_{\hat{W}}(v))$
where $\mu_{\hat{W}}:\hat{W}_{p}\to \lih^*$ is the moment map for the linear $Stab\, p$ action and $\lih=Lie(Stab\, p)$.
\end{theorem}
As a consequence, a neighborhood of $x\in \mu^{-1}(0)/G$,
corresponding to $p$, is homeomorphic to $\hat{W}_{p}\git Stab\,p$
which is the reduction of $G\times_{Stab\,p}((\lig/\lih)^*\times\hat{W}_{p})$
by $G$.

In the case of extended moduli space, we showed in the previous section
that the symplectic slice at a point 
$\phi\in \mathrm{Hom}(\pi,G)$
is $\hat{W}_{\phi}=H^1(\pi,\lig_{Ad\phi})$
on which $Stab\, p$ acts by conjugation 
on $\lig$. The moment map for this action
 is given by $\mu_{\hat{W}_{\phi}}(u)=[u,u]$, 
$\mu_{\hat{W}_{\phi}}:H^1(\pi,\lig_{Ad\phi})
\to H^2(\pi,\lig_{Ad\phi})\cong \lih^*$ using the Lie bracket of $\lig$
and the cup product. 

The above local normal form theorem now shows that 
in a neighborhood of the image of a point $\phi\in \mathrm{Hom}(\pi,G)$,
$X(G)$ is homeomorphic to the symplectic reduction 
$H^1(\pi,\lig_{Ad\phi})\git Stab\, \phi$. This is again consistent with
the deformation theory result of Goldman-Millson \cite{GolMil} 
but our approach is more in the spirit of
symplectic geometry.

Let $\lih=Lie(Stab\,\phi)$ and $\mathfrak{z}(\lih)=
\{\zeta\,|\,[\zeta,\zeta']=0 \text{ for all }\zeta'\in \lih\}$.
Recall that $H^1(\pi,\lig_{Ad\phi})$ is the quotient of 
the space of cocycles
\begin{equation}\label{cycle}
Z^1(\pi,\lig_{Ad\phi})=\{u:\pi\to\lig\,|\,u(rs)=u(r)+Ad \phi (r) u(s)\}
\end{equation}
by the space of coboundaries
\begin{equation}
B^1(\pi,\lig_{Ad\phi})=\{u:\pi\to\lig\,|\,u(r)=Ad \phi(r) u_0-u_0
\text{ for } u_0\in \lig\}.\end{equation}

Since $Ad\phi$ stabilizes $\lih$, $Ad \phi$ preserves 
$\mathfrak{z}(\lih)$ and $\mathfrak{z}(\lih)^{\perp}$ for 
the orthogonal decomposition $\lig=\mathfrak{z}(\lih)
\oplus \mathfrak{z}(\lih)^{\perp}$. 
Hence, each $u:\pi\to\lig$ is the sum
of cocycles $u_{\lih}:\pi\to\lig\to\mathfrak{z}(\lih)$ and $u_{\lih}^{\perp}:
\pi\to\lig\to \mathfrak{z}(\lih)^{\perp}$.
Thus, we get the induced decomposition
$$H^1(\pi,\lig_{Ad \phi})= H^1(\pi,\mathfrak{z}(\lih)_{Ad\phi})
\oplus H^1(\pi,\mathfrak{z}(\lih)^{\perp}_{Ad\phi}).$$ 
In terms of gauge theory, this is obtained by 
reduction of the structure group to the centralizer $H'$ of
the identity component $H$ of $Stab(\phi)$.
Clearly, the first factor is $\lih$ fixed while the second factor
has no nonzero fixed point.\footnote{One easy way to see it is
to consider $H^1(\pi,\mathfrak{z}(\lih)^{\perp}_{Ad\phi})$ as
the space of harmonic $\mathfrak{z}(\lih)^{\perp}$ valued 1-forms.}
 Hence, the $\lih$-fixed subspace
in $\hat{W}_{\phi}$ is 
$\hat{W}^{\lih}_{\phi}:=H^1(\pi,\mathfrak{z}(\lih)_{Ad\phi})$
and its orthogonal complement is 
$W_{\phi}:=H^1(\pi,\mathfrak{z}(\lih)^{\perp}_{Ad\phi})$

In fact, the (orbifold) stratum direction 
$\hat{W}^{\lih}_{\phi}$ can be integrated: 
since $\lih =\mathfrak{z}(\phi)$ and $\mathfrak{z}(\mathfrak{z}(\lih))=\lih$,
 $\phi$ is an irreducible point\footnote{Here, ``irreducible''
means that the infinitesimal stabilizer is of minimal dimension.}
 in $\mathrm{Hom}(\pi,H')$
where $H'$ is the centralizer of $H=Stab(\phi)$. We get 
the following decompositions 
$$\mathrm{Hom}(\pi,G)=\cup_{(H)\in I}G \mathrm{Hom}(\pi,H')^{irr}$$
$$\mathrm{Hom}(\pi,G)/G=\cup_{(H)\in I} \mathrm{Hom}(\pi,H')^{irr}/N^{H'}$$
whose $irr$ denotes the irreducible part.

\section{Weakly balanced action and the splitting theorem}
In this section, we show that the $G$ action on the extended moduli
space is weakly balanced in the sense of \cite{kiem,kw} and
hence we get an isomorphism of the truncated equivariant cohomology
of $\mathrm{Hom}(\pi,G)$ with the intersection cohomology of 
the representation space $X(G)$.

First, we recall the weakly balanced condition. 
Let $U$ be a Hamiltonian $G$ space with moment map $\mu$ and $Z=\mu^{-1}(0)$.
Let $H$ be the identity component of $Stab\,p$ for a point $p\in Z$ and 
$\hat{W}_p$ denote the symplectic slice at $p$.
We decompose $\hat{W}_p=\hat{W}_p^{\lih}\oplus W_p$
into the $H$ fixed subspace and its complement.
 Choose a compatible 
complex structure on the symplectic vector space $W_{p}$
 and let $\pp W_{p}$ denote the complex projective space associated
to $W_{p}$. The unstable strata in $\pp W_{p}$ are parametrized
by a finite set $\mathcal{B}$ in the Weyl chamber. Namely, $\mathcal{B}$
is the set of the closest points from the origin to the convex hulls
of some weights of the action of $H$ on $W_{p}$. (See \cite{K1}.)
Given $\beta\in\mathcal{B}$, let $n(\beta)$ denote the number of 
weights $\alpha$ such that $\langle\alpha,\beta\rangle\,<\,
\langle\beta,\beta\rangle$
with respect to the invariant metric.

\begin{definition}
We say the $H$ action on $W_{p}$ is \emph{linearly balanced} if 
$2n(\beta)\ge \frac12 \dim W_{p}$. 
\end{definition}
This condition makes the codimensions of the unstable strata
large enough. 

For any subgroup $K\subset H$ which also appears as the identity component 
of the stabilizer of a point in $Z$, 
we consider the $N^{K}\cap H/K_{}$ action
on the $K$ fixed point subspace $W_{p}^{\lik}$ where  
$N^{K}$ is the normalizer of $K$.
\begin{definition}
We say the $G$ action on a Hamiltonian space $U$ is \emph{balanced}
if for all $p\in Z$ the action of the identity
component $H$ of its stabilizer on $W_{p}$ is linearly balanced
and if for a subgroup $K$ as above the $N^{K}\cap H/{K}$ action
on $W_{p}^{\lik}$ is also linearly balanced.
\end{definition}

The usefulness of this condition becomes evident from the splitting theorem
in \cite{kiem,kw}. To describe this theorem, we need to recall
the ``truncated equivariant cohomology''. For each 
connected subgroup $H$ of $G$, which is the identity component 
of the stabilizer of a point in $Z$, let $Y_{H}$ denote
the subset of $Z$ fixed by $H$.
Consider the obvious map 
$$GY_{H}\leftarrow G\times_{N^{H}}Y_{H}.$$
This composed with the inclusion induces a map of equivariant
cohomology groups
$$H^*_G(Z)\to H^*_G(GY_{H})\to H^*_{N^{H}}
(Y_{H})$$ which we denote by $\Psi_{H}$.
Note that $H^*_{N^{H}}(Y_{H})\cong
[H^*_{N^{H}_0/H}(Y_{H})\otimes H^*_{H}]^{\pi_0N^{H}}$
where $H^*_{H}$ denotes the cohomology of the classifying space of $H$
and $N^{H}_0$ is the identity component of $N^{H}$.
We define the truncated equivariant cohomology  by
$$V^*_U:=\{\xi\in H^*_G(Z)\,|\, \Psi_{H}(\xi)\in 
H^*_{N^{H}_0/H}(Y_{H})\otimes H^{<n_{H}}_{H}\}$$
where $n_{H}=\frac12 \dim (W_{p}\git H)$
for $p\in Z$ whose stabilizer has $H$ as its identity component.

\begin{theorem}\cite{kiem, kw}
Suppose the proper Hamiltonian $G$ space $U$ is balanced. Then there is
an isomorphism $\rho:IH^*(X)\to V^*_U$ such that the
intersection pairing $\langle,\rangle$ of $IH^*(X)$ is given by 
$$\rho(\alpha)\cup \rho(\beta)=\langle \alpha, \beta\rangle \rho(\tau)$$
for all $\alpha, \beta \in IH^*(X)$ of complementary degrees
with respect to the real dimension of $X=U\git G$ where $\tau$ is the
fundamental class.\end{theorem}

In fact, the theorem is proved in a more general context, namely
under the ``weakly balanced'' assumption. For our purpose in this
paper, we will need only the balanced condition.

Let us now return to  the extended moduli space. For each point $\phi$
in $\mu^{-1}(0)=\mathrm{Hom}(\pi,G)$, we have $W_{\phi}=H^1(\pi,
\mathfrak{z}(\lih)^{\perp}_{Ad\phi})$ where $Stab\,\phi$ acts on 
$\mathfrak{z}(\lih)^{\perp}$ by conjugation and $\lih=Lie(Stab\,\phi)$.
There is a geometric way to assign a compatible complex structure
on this space. Choose a complex structure on $\Sigma$
and consider the holomorphic  vector bundle over $\Sigma$
with fiber $\mathfrak{z}(\lih)^{\perp}\otimes\cc$ given by $\phi$
$$E=\tilde{\Sigma}\times_{\pi}(\mathfrak{z}(\lih)^{\perp}\otimes\cc)
\to \Sigma$$ where $\tilde{\Sigma}$ is the universal cover.
Then $W_{\phi}$ is canonically isomorphic to the complex vector space
 $H^{0,1}(\Sigma,\mathfrak{z}(\lih)^{\perp}\otimes\cc)$.

Consider the root space decomposition of the complex semisimple
Lie algebra $\lig\otimes \cc$ with respect to a Cartan
subalgebra of $\lig$ containing a Cartan subalgebra $\lit_{\lih}$ of $\lih$. 
Of course, $\lambda$ is a root iff $-\lambda$ is a root. Hence,
the weights of the $\lit_{\lih}$ action on $\lig$ are
symmetric with respect to the origin and clearly 
$\mathfrak{z}(\lih)\otimes\cc$ lies in the zero eigenspace.
Therefore, the weights of the $\lit_{\lih}$ action on 
$\mathfrak{z}(\lih)^{\perp}\otimes\cc$ are symmetric with
respect to the origin. Now, because $\lit_{\lih}$ commutes
with $\phi$, the eigenspace decomposition of the fiber 
$\mathfrak{z}(\lih)^{\perp}\otimes\cc=\bigoplus_{\lambda} F_{\lambda}$ 
provides us with a decomposition $E=\bigoplus_{\lambda} E_{\lambda}$ 
of the vector bundle into subbundles of eigenspaces. 
The degrees of the $E_{\lambda}$ are all zero and 
by the invariant metric of $\lig$, $E_{\lambda}\cong E_{-\lambda}^*$. 
In particular, the subbundles $E_{\lambda}, E_{-\lambda}$ have the same
Riemann-Roch number. Moreover,  
$H^0(\Sigma,E_{\lambda})\cong H^0(\pi,F_{\lambda})=\{v\in F_{\lambda}\,|\,
Ad\phi \text{ fixes } v\}$ gives rise to the maximal
 trivial subbundle of $E_{\lambda}$
and a similar description for $E_{-\lambda}$ gives us the
equality $\dim H^0(\Sigma,E_{\lambda})=\dim H^0(\Sigma,E_{-\lambda})$
via the isomorphism $E_{-\lambda}\cong E^*_{\lambda}$.
Therefore, we conclude that 
$\dim H^{0,1}(\Sigma,E_{\lambda})=\dim H^{0,1}(\Sigma,E_{-\lambda})$ and thus 
the weights of the $\lit_{\lih}$ action on $W_{\lih}$ are symmetric
with respect to the origin. In particular, the action is 
linearly balanced.

Next, we consider the $\lik$ fixed point subspace $W_{\phi}^{\lik}$ for
 another infinitesimal stabilizer $\lik \subset\lih$. Let $K$
be the connected subgroup of $Stab\, \phi$ whose Lie algebra is $\lik$. 
Notice that in our case 
$$W_{\phi}^{\lik}=H^1(\pi,(\mathfrak{z}(\lih)^{\perp}
\cap\mathfrak{z}(\lik))_{Ad\phi}).$$
$\mathfrak{z}(\lik)\otimes \cc$ is an ideal in the
complexification of the Lie algebra of $N^{K}$
and thus the weights of the $N^{K}$ action on $\mathfrak{z}(\lik)\otimes\cc$
by conjugation are symmetric with respect to the origin. 
As $K$ fixes $\mathfrak{z}(\lik)\otimes\cc$, the same is true for the
action of $N^{K}/K$ on $\mathfrak{z}(\lik)\otimes \cc$
and hence so is the action of the subgroup $H\cap N^{K}/K$.
Now, $H\cap N^{K}/K$ acts trivially on $\mathfrak{z}(\lih)\otimes\cc$
and therefore the weights of
the $N^{K}\cap H/{K}$ action on $(\mathfrak{z}(\lih)^{\perp}
\cap\mathfrak{z}(\lik))\otimes \cc$ are symmetric with respect 
to the origin. As in the previous paragraph, we can 
deduce that the $N^{K}\cap H/{K}$ action on $W_{\phi}^{\lik}$
is also linearly balanced.
In summary, we proved the following theorem.

\begin{theorem}
The $G$ action on $U=\mu^{-1}(B)\subset \mathcal{N}^{\lig}$
is balanced. Therefore, we have the isomorphism
$$IH^*(X(G))\cong V^*_{U}\subset H^*_G(\mathrm{Hom}(\pi,G)).$$
\end{theorem}

Moreover, the above isomorphism preserves the mapping class 
group action since it is defined by pulling back 
certain differential forms \cite{kiem} via the quotient map $U\to X$ which
is equivariant with respect to the mapping class group 
$\Gamma_{g,1}$ of the punctured Riemann surface.

The equivariant cohomology $H^*_G(\mathrm{Hom}(\pi,G))$
is isomorphic to the gauge group equivariant cohomology
of the space of flat $G$ connections on $\Sigma$,
which can be computed by the Morse theory of \cite{ab2}. 
The cohomology ring $H^*_{SU(2)}(\mathrm{Hom}
(\pi,SU(2)))$ was completely determined in \cite{kiem4} and 
hence the above theorem can be applied to compute
the intersection cohomology of the representation space
as we will see in the next section.

\section{$SU(2)$ case}
In this section, we focus on the case where $G=SU(2)$.
First, we recall the structure theorem of the equivariant
cohomology which determines the cup product structure.

Let $\tilde \Sigma$ be the universal cover of $\Sigma$. Consider
$$\tilde{\Sigma}\times_{\pi}(\mathrm{Hom}(\pi,SU(2))\times End(\Bbb C^2))
\rightarrow \Sigma\times \mathrm{Hom}(\pi,SU(2))$$
where $g\in \pi$ maps $(\phi, v)\in \mathrm{Hom}(\pi,SU(2))
\times End(\Bbb C^2)$
to $(\phi, Ad\,\phi(g)\,v).$ It is a vector bundle of rank 4 which induces
$End\,\Cal U$ over 
$\Sigma\times (\mathrm{Hom}(\pi,SU(2))\times_{PU(2)}EPU(2))$ by pulling back
and taking quotient.
We define equivariant classes $\alpha, \beta, \psi_i$ of degree 2,4,3
respectively by
$$c_2(End(\Cal U))=2\alpha \otimes [\Sigma]+4\sum_{i=1}^{2g} \psi_i \otimes 
e_i -\beta \otimes 1$$
where $e_i$ is a symplectic basis of $H^1(\Sigma)$ so that $e_ie_{i+g}=
[\Sigma]$.

\begin{theorem} \cite{kiem4} \label{struthm}
Let $\gamma=-2\sum_{i=1}^g
\psi_i\psi_{i+g}$ and 
consider the Lefschetz decompositon of the exterior algebra
$$\wedge(\psi_1,\cdots, \psi_{2g})=\bigoplus_{l=1}^g Prim_l\otimes
\cc[\gamma]/\gamma^{g-l+1}$$
where $Prim_l$ is the degree $3l$ primitive part.
Then we have 
$$H^*_{SU(2)}(\mathrm{Hom}(\pi,SU(2)))\cong
\oplus_{l=0}^gPrim_l\otimes \cc [\alpha,\beta,\gamma]/I_{g-l}.$$
Here $I_k$ is the ideal of the 
polynomial ring $\cc[\alpha,\beta,\gamma]$ generated by
$c_{k+1},c_{k+2},c_{k+3}$ where the sequence $\{c_n\}$ is defined
by $$nc_n=\alpha c_{n-1}+(n-2)\beta c_{n-2}+2\gamma c_{n-3}$$
with $c_0=1,\, c_1=\alpha,\, c_2=\frac{\alpha^2}{2}$, etc.
\end{theorem}
 
Because $Prim_l$s are the irreducible $Sp_{2g}$ modules, 
this theorem completely describes the mapping class group action
as it factors through the symplectic group action on $\psi_i$. 
Furthermore, a Gr\"obner basis for $I_k$ was given in \cite{kiem4} so that
the cup product can be computed efficiently. 

Now by the splitting theorem
above, $IH^*(X(SU(2)))$ is isomorphic to $V^*_U$ which we intend to compute.
Note that up to conjugation 
$\{SU(2), S^1\}$ are the only possible stabilizers of
points in $\mathrm{Hom}(\pi,SU(2))$.
It turns out that we don't have to think about truncation on the $SU(2)$ 
fixed point set because it is taken care of by that for $S^1$. For that
one has only to observe 
that $Y_{SU(2)}\subset Y_{S^1}$ and that 
$n_{SU(2)}=3g-3>n_{S^1}=2g-3$.

We have to consider the following map
$$H^*_{SU(2)}(\mathrm{Hom}(\pi,SU(2)))\rightarrow 
H^*_{SU(2)}(SU(2)\times_{N^{S^1}}
Y_{S^1})=[H^*(Jac)\otimes \cc [u]]^{\Bbb Z/2}$$
where $Jac$ denotes the Jacobian variety and $\Bbb Z/2$ acts as $-1$ on both
components. Let $d_i$ be the basis of $H^1(Jac)$ defined as the K\"unneth
coefficients of the first Chern class of a universal line bundle corresponding
to $e_i$. 
According to \cite{clm}, $\alpha$ is mapped to 
$w=-2\sum_{i=1}^gd_id_{i+g}$, $\beta$ to $4u^2$, 
and $\psi_i$ to $-2u d_i$. \par

Let $P_t^{SU(2)}(\mathrm{Hom}(\pi,SU(2)))$ denote
the Poicar\'e series for the equivariant cohomology $H^*_{SU(2)}(
\mathrm{Hom}(\pi,SU(2)))$. Then it is well known from \cite{ab2} that 
$$P_t^{SU(2)}(\mathrm{Hom}(\pi,SU(2)))=
\frac{(1+t^3)^{2g}-t^{2g+2}(1+t)^{2g}}{(1-t^2)(1-t^4)}.$$
For the intersection Poicar\'e series $IP_t(X)=\sum t^i\dim IH^i(X)$, we have only to subtract out the Poincare series of 
$$Im\{H^*_{SU(2)}(\mathrm{Hom}(\pi,SU(2)))\rightarrow 
[H^*(Jac)\otimes \Bbb C[u]]^{\Bbb Z/2}\}\,\,\,\cap\,\,\, 
\{H^*(Jac)\otimes u^{g-1}\Bbb C[u]\}.$$
By the Lefshetz decomposition and a trivial combinatorial argument, 
one can easily see that the intersection is actually
$$[H^*(Jac)\otimes u^{g-1}\Bbb C[u]]^{\Bbb Z/2}$$
whose Poincare series is precisely
$$\frac12\{\frac{(1+t)^{2g}(t^2)^{g-1}}{1-t^2}+
\frac{(1-t)^{2g}(-t^2)^{g-1}}{1+t^2}\}.$$

\begin{proposition}\label{bettisu} $$IP_t(X(SU(2)))=
P_t^{SU(2)}(\mathrm{Hom}(\pi,SU(2)))-
\frac12\{\frac{(1+t)^{2g}(t^2)^{g-1}}{1-t^2}+
\frac{(1-t)^{2g}(-t^2)^{g-1}}{1+t^2}\}.$$\end{proposition}

It is an elementary exercise to check that Proposition \ref{bettisu}
indeed coincides with Kirwan's computations \cite{k5}.

A little bit more careful examination of the truncation map
together with Theorem \ref{struthm} described above gives us the following.
\begin{theorem} Let $E_m$ be the vector space spanned by
$$\{\alpha^i\beta^j\xi^k\,|\,
\text{(1) }i+2k\le m, \text{ (2) }j+2k\le m,
\text{ (3) if }k=0\text{ then }j<[\frac{m}{2}]\}$$
where $\xi=\alpha\beta+2\gamma$.
Then $$IH^*(X(SU(2)))\cong \oplus_{i=0}^g Prim_l\otimes E_{g-l}.$$
\end{theorem}

Again, this theorem precisely describes the action of the mapping class group
$\Gamma_{g,1}$ 
as our isomorphism is equivariant and hence the action also factors through the symplectic group action on $Prim_l$.

Also the intersection pairing can be computed in terms of the
cup product structure. One can prove that 
$$\alpha^m\beta^n=-m!b_{g-n-1}\frac{\alpha^{g-2}
\beta^{g-2}\xi}{(g-2)!}$$
where $m+2n=3g-3$, $n<g-1$ and $b_k$ are given by
$$\frac{t}{\mathrm{tanh} t}=\sum_{k\ge 0}b_kt^{2k},$$
directly from the structure theorem at least
for low genus case.
Hence, if we take $\frac{\alpha^{g-2}\beta^{g-2}\xi}{(g-2)!(-4)^{g-1}}$,
as our fundamental class then  we get 
$$\langle \kappa(\alpha^i\beta^j),\kappa(\alpha^k\beta^l)
\rangle=-(-4)^{g-1}m!b_{g-n-1}$$
for $i+k=m$, $j+l=n$, $m+2n=3g-3$, $n<g-1$.
In principle, it is a number theoretic or combinatorial exercise
to deduce the above formula from the structure theorem but
it seems very difficult to achieve in practice. 

However the above pairing formula can be justified by using 
the computation of intersection numbers
on the moduli spaces of parabolic bundles by Jeffrey and Kirwan
\cite{JK2}. This will be done in \cite{JKKW2}.

\end{document}